\documentclass{amsart}
\usepackage{graphicx}
\theoremstyle{definition}

\theoremstyle{remark}

\numberwithin{equation}{section}
\begin{document}

\title{A note on the Ricci flow on noncompact manifolds}%
\author{Hong Huang}%
\address{School of Mathematical Sciences,Key Laboratory of Mathematics and Complex Systems,Beijing Normal University,Beijing 100875, P. R. China}%
\email{hhuang@bnu.edu.cn}%

\thanks{Partially supported by NSFC no.10671018.}%


%\date{}%
%\dedicatory{}%
%\commby{}%
% ----------------------------------------------------------------
\begin{abstract}
Let $(M^3,g_0)$ be a complete noncompact Riemannian 3-manifold
with nonnegative Ricci curvature and with injectivity radius
bounded away from zero. Suppose that the scalar curvature
$R(x)\rightarrow 0$ as $x\rightarrow \infty$. Then the Ricci flow
with initial data $(M^3,g_0)$ has a long time solution on $M^3
\times [0,\infty)$. This extends a recent result of Ma and Zhu. We
also have a higher dimensional version, and we reprove a
K$\ddot{a}$hler analogy due to Chau, Tam and Yu.
\end{abstract} \maketitle
% ----------------------------------------------------------------

\section {Introduction}

In a recent paper [MZ] Ma and Zhu announced  the following result:
Let $(M^3,g_0)$ be a complete noncompact Riemannian 3-manifold
with nonnegative sectional curvature such that $|Rm|(x)\rightarrow
0$ as $x\rightarrow \infty$. Then the Ricci flow with initial data
$(M^3,g_0)$ has a long time solution on $M^3 \times [0,\infty)$.

In this short note we'll improve Ma and Zhu's result by showing the
following

\hspace *{0.4cm}

 {\bf Theorem 1} Let $(M^3,g_0)$ be a complete
noncompact Riemannian 3-manifold with nonnegative Ricci curvature
and with injectivity radius bounded away from zero. Further assume
that the scalar curvature $R(x)\rightarrow 0$ as $x\rightarrow
\infty$. Then the Ricci flow with initial data $(M^3,g_0)$ has a
long time solution on $M^3 \times [0,\infty)$.

\hspace *{0.4cm}

 Similarly one has

\hspace *{0.4cm}

 {\bf Theorem 2} Let $(M^n,g_0)$ be a complete
noncompact Riemannian $n$-manifold with nonnegative curvature
operator and with injectivity radius bounded away from zero.
Further assume
 that the scalar curvature $R(x)\rightarrow 0$ as $x\rightarrow
\infty$. Then the Ricci flow with initial data $(M^n,g_0)$ has a
long time solution on $M^n \times [0,\infty)$.

(The case $n=3$ is due to Ma and Zhu ([MZ, Theorem 1.1]),  where
the injectivity radius condition is not assumed, but actually this
condition does not follow from the nonnegativity of curvature
operator.)

 \hspace *{0.4cm}

We also reprove a K$\ddot{a}$hler analogy due to Chau, Tam and Yu.

\hspace *{0.4cm}

 {\bf Theorem 3} ([CTY,Theorem 1.2]) Let $(M^n, g_0)$ be a complete noncompact
 K$\ddot{a}$hler manifold with nonnegative holomorphic bisectional
 curvature and with injectivity radius bounded away from zero. Further
assume that  $|Rm|(x)\rightarrow 0$ as $x\rightarrow \infty$. Then
the K$\ddot{a}$hler-Ricci flow with initial data $(M^n,g_0)$ has a
long time solution on $M^n \times [0,\infty)$.

(Actually the condition $|Rm|(x)\rightarrow 0$ as $x\rightarrow
\infty$ may be replaced by $R(x)\rightarrow 0$ as $x\rightarrow
\infty$.)

 \hspace *{0.4cm}

Both theorems are consequences of Perelman's no local collapsing
theorem, pseudolocality theorem (actually  its generalization due
to Chau, Tam and Yu [CTY]),  and his theorem on asymptotic volume
ratio on ancient solution to Ricci flow ( or its K$\ddot{a}$hler
analogy due to Ni [N]), as shown in the next section.

\section {Proof of Theorems}

 Proof of Theorem 1.
 First note that since $(M^3,g_0)$  has nonnegative Ricci curvature and
 the scalar  curvature $R(x)\rightarrow 0$  as $x\rightarrow
\infty$, we have $|Rm|(x)\rightarrow 0$ as $x\rightarrow \infty$, so
the sectional curvature of $(M^3,g_0)$ is bounded, then we have the
short time existence of a solution $(M^3,g(t))$ to the Ricci flow
equation $\frac{\partial g_{ij}}{\partial t}=-2R_{ij}$ with initial
data $(M^3,g_0)$ by Shi([S2]). Also note that the Ricci curvature of
$g(t)$ remains nonnegative along the Ricci flow again by Shi ([S1]).

Let the maximal time interval of existence of $(M^3,g(t))$ be
$[0,T)$. We'll prove by contradiction that $T=\infty$. Suppose
$T<\infty$.  Then since $(M^3,g_0)$ has bounded sectional curvature
and has injectivity radius bounded away from zero, by the virtue of
Perelman's work [P], see Kleiner-Lott [KL,Theorem 26.2],
$(M^3,g(t))$ is non-collapsed (in the sense of [KL,Definition
26.1]).

Since $(M^3,g_0)$ has injectivity radius bounded away from zero, and
satisfies $|Rm|(x)\rightarrow 0$ as $x\rightarrow \infty$, and $T<
\infty$, by a theorem of Chau, Tam and Yu([CTY, Theorem 1.1] which
extends Perelman's pseudolocality theorem in [P]) we know that there
exists some compact set $S \subset M^3$ such that $|Rm|(x,t)$ is
uniformly bounded on $({M^{3}}\setminus S) \times [0,T)$. But the
sectional curvature of $(M^3, g(t))$ must blow up at time $T$, hence
there exist $t_n \rightarrow T$, $p_n \in S$ such that
$Q_n:=|Rm|(p_n,t_n)=sup_{x \in M^3, t\leq t_n} |Rm|(x,t)\rightarrow
\infty$ as $n\rightarrow \infty$.  By Hamilton's compactness theorem
for Ricci flow ([H]), the rescaled solutions $(M^3,Q_n
g(t_n+{Q_n}^{-1}t),p_n)$ sub-converge to a complete non-flat ancient
solution $(M_\infty,g_\infty(t),q)$ to  Ricci flow, which has
(bounded) nonnegative sectional curvature by the Hamilton-Ivey
curvature pinching estimate. Then $(M_\infty,g_\infty(t))$ has zero
asymptotic volume ratio, again by Perelman [P]. So given any
$\varepsilon
>0$, there exists $r>0$ such that

$\frac{Vol B_{g_{\infty}(0)}(q,r)}{r^3}<\frac{\varepsilon}{2}$.

Then, when $n$ is sufficiently large we have

$\frac{Vol B_{g(t_n)}(p_n, {Q_n}^{-\frac{1}{2}}r)}
{({Q_n}^{-\frac{1}{2}}r)^3}< \varepsilon$.     \  \  \  \ \ \ (*)

Since $|Ric|(x,t)$ is uniformly bounded on $(M^3 \setminus S)\times
[0,T)$, the volume of any compact region in $M^3\setminus S$ decays
with a controllable speed. So there exist a positive constant
$\delta$ and a compact region $\Omega \subset M^3\setminus S$ such
that $Vol_{g(t)}\Omega \geq \delta$, $t \in [0,T)$.

Choose $r_0$ sufficiently large such that $\Omega \subset
B_{g_0}(p,r_0)$ for any $p \in S$. Then since the Ricci curvature of
$g(t)$ remains nonnegative, we know that the distance function of
$(M^3, g(t))$ is decreasing, so we have $\Omega \subset
B_{g(t)}(p,r_0)$, for any $t \in [0,T)$, and any $p \in S$.

Choose $n$ sufficiently large such that ${Q_n}^{-\frac{1}{2}}r<r_0$.
Then by Bishop-Gromov theorem we have

$\frac{Vol B_{g(t_n)}(p_n, {Q_n}^{-\frac{1}{2}}r)}
{({Q_n}^{-\frac{1}{2}}r)^3} \geq \frac{Vol B_{g(t_n)}(p_n,r_0)}
{{r_0}^3}\geq \frac{\delta}{{r_0}^3}$.

But if we choose $\varepsilon <\delta/{r_0}^3$, then the above
inequality contradicts the inequality (*), and we are done.

\hspace *{0.4cm}

{\bf Remark } \ \ The above proof is a  modification of that of Ma
and Zhu ([MZ]). It seems to me that the arguments in [MZ] may have
some gaps. In their arguments in [MZ], Ma and Zhu invoke  a theorem
of Hamilton (see Theorem 2.3 in [MZ]) to imply the existence of $S$
( where it is denoted by $K$), but it seems that Hamilton's theorem
is not sufficient for this implication.

\hspace *{0.4cm}

The prove of Theorem 2 is similar, and is omitted. (When the
condition in Theorem 2 on nonnegativity of curvature operator is
replaced by positivity of curvature operator, then the injectivity
radius condition may be removed, and the same result holds, since
in this case one has the  Gromoll-Meyer injectivity radius
estimate. Actually in this case the proof is simpler, since one
does not need Perelman's no local collapsing theorem.)

\hspace *{0.4cm}

To prove Theorem 3, one only needs to replace Perelman's theorem
on asymptotic volume ratio by Ni's([N,Theorem 2]), which says that
a non-flat ancient solution to K$\ddot{a}$hler-Ricci flow with
bounded nonnegative holomorphic bisectional curvature must have
zero asymptotic volume ratio. Then after minor modifications the
arguments in the proof of Theorem 1 apply to this case.

% ----------------------------------------------------------------
\bibliographystyle{amsplain}

\hspace *{0.4cm}

{\bf Reference}

\bibliography{1}[CTY] A. Chau, L.-F. Tam, and C. Yu, Pseudolocality for
the Ricci flow and applications, arXiv:math/0701153v2.

\bibliography{2}[H] R. Hamilton, A compactness property for
solutions of the Ricci flow, Amer. J. Math. 117 (1995), no.3,
545-572.

\bibliography{3}[KL] B. Kleiner and J. Lott, Notes on Perelman's
papers, arXiv:math/0605667v2.

\bibliography{4}[MZ] L. Ma and A. Zhu, Nonsingular Ricci flow on a noncompact manifold
in dimension three, arXiv:0806.4672v1.

\bibliography{5}[N] L. Ni, Ancient solutions to
K$\ddot{a}$hler-Ricci flow, Math. Res. Lett. 12 (2005), no. 5-6,
633-653.

\bibliography{6}[P] G. Perelman, The entropy formula for the Ricci
flow and its geometric applications, arXiv:math.DG/0211159v1.

\bibliography{7}[S1]W.X. Shi, Complete noncompact three manifolds with
nonnegative Ricci curvature, J. Diff. Geom. 29(1989),353-360.

\bibliography{8}[S2] W.X. Shi, Deforming the metric on complete Riemannian manifolds,
J. Diff. Geom. 30 (1989),no.1, 223-301.

\end{document}